\documentclass{conm-p-l}
\usepackage{amssymb,amsxtra}

\newcommand{\norm}[1]{\mbox{$\|#1\|$}}

\newcommand{\x}{\times}
\newcommand{\cs}{\mbox{$C^{*}$-algebra}}
\newcommand{\css}{\mbox{$C^{*}$-algebras}}

\newcommand{\C}{\mathbb{C}}

\newcommand{\R}{\mathbb{R}}

\newcommand{\ov}[1]{\mbox{$\overline{#1}$}}

\newcommand{\al}{\mbox{$\alpha$}}

\newcommand{\bt}{\mbox{$\beta$}}

\newcommand{\Ga}{\mbox{$\Gamma$}}

\newcommand{\De}{\mbox{$\Delta$}}

\newcommand{\la}{\mbox{$\lambda$}}

\newcommand{\mfB}{\mathfrak{B}}

\newcommand{\mfK}{\mathfrak{K}}
\newcommand{\mfL}{\mathfrak{L}}

\newcommand{\mfH}{\mathfrak{H}}

\newcommand{\bc}{\begin{center}}

\newcommand{\ec}{\end{center}}
\newcommand{\be}{\begin{enumerate}}
\newcommand{\ee}{\end{enumerate}}

\newcommand{\beqn}{\begin{eqnarray}}
\newcommand{\eeqn}{\end{eqnarray}}
\newcommand{\beqns}{\begin{eqnarray*}}
\newcommand{\eeqns}{\end{eqnarray*}}
\newcommand{\bq}{\begin{quote}}
\newcommand{\eq}{\end{quote}}

\newcommand{\bi}{\begin{itemize}}
\newcommand{\ei}{\end{itemize}}
\newcommand{\bd}{\begin{description}}
\newcommand{\ed}{\end{description}}

\newcommand{\lan}{\mbox{$\langle$}}
\newcommand{\ran}{\mbox{$\rangle$}}

\theoremstyle{plain}

\theoremstyle{definition}

\theoremstyle{remark}

\numberwithin{equation}{section}

\begin{document}

\title{The Fourier-Stieltjes and Fourier algebras for locally compact groupoids}

\author{Alan L. T. Paterson}
\address{Department of Mathematics, University of Mississippi, University,
Mississippi 38677}
\email{mmap@olemiss.edu}
\date{January 31, 2002.}

\begin{abstract}
The Fourier-Stieltjes and Fourier algebras $B(G)$, $A(G)$ for a general locally compact group $G$, first studied by P. Eymard,  have played an important role in harmonic analysis and in the study of the operator algebras generated by $G$.  Recently, there has been interest in developing versions of these algebras for locally compact {\em groupoids}, justification for this being that, just as in the group case, the algebras should play a useful role in the study of groupoid operator algebras.  Versions of these algebras for the locally compact groupoid case appear in three related theories: (1) a measured groupoid theory (J. Renault), (2) a Borel theory (A. Ramsay and M. Walter), and (3) a continuous theory (A. Paterson).  
The present paper is expository in character.  For motivational reasons, it starts with a description of the theory of $B(G)$, $A(G)$ in the locally compact group case, before discussing these three theories.  Some open questions are also raised.
\end{abstract}
\maketitle

\section{$A(G)$ and $B(G)$ for locally compact groups}

I am grateful to the organizers of the 2001 Conference on Banach spaces at the University of Memphis for the opportunity to participate in that exciting conference.  

We first specify some notation.  If $X$ is a locally compact Hausdorff space, then $C(X)$ is the algebra of bounded, continuous, complex-valued functions on $X$.  The space of functions in $C(X)$ that vanish at $\infty$ is denoted by $C_{0}(X)$, while $C_{c}(X)$ is the space of functions in $C(X)$ with compact support.  The space of complex, bounded, regular Borel measures on $X$ is denoted by $M(X)$.

The Banach spaces $B(G)$, $A(G)$ (where $G$ is a locally compact groupoid) that we will consider in this paper arise naturally in the group case in non-commutative harmonic analysis and duality theory.  (See later in this section.)  When $G$ is a locally compact group, $B(G)$ and $A(G)$ are just the (much studied) Fourier-Stieltjes and Fourier algebras.  The need to have available versions of these Banach algebras for the case of a locally compact groupoid is a consequence of the fact that many of the operator algebras of present day interest, especially in non-commutative geometry, come from  groupoid, and not (in any obvious way) from group representations, so that having available versions of $B(G)$, $A(G)$ in the groupoid case would provide the resources for extending to that case
the properties of group operator algebras that depend on $B(G)$, $A(G)$.

We start by discussing these algebras in the locally compact group case.
Let us, for the present, specialize even further by letting $G$ be a locally compact {\em abelian} group with character space $\widehat{G}$.  An element of 
$\widehat{G}$ is a continuous homomorphism $t:G\to {\bf T}$, and $\widehat{G}$ is a locally compact abelian group with pointwise product and the topology of uniform convergence on compacta. The Fourier transform 
$f\to \hat{f}$ takes $f\in L^{1}(G)$ into $C_{0}(\widehat{G})$, where 
\[  \hat{f}(t)=\int f(x)\ov{t(x)}\,dx  \]
where $dx$ is a left Haar measure on $G$.  The inverse Fourier transform
$\mu\to \check{\mu}$ takes $M(\widehat{G})$ back into $C(G)$ where 
\[  \check{\mu}(x)=\int \hat{x}(t)\,d\mu(t).            \]
For example, when $G=\R$, we have $\widehat{G}=\R$ where $t\in \widehat{G}$ is associated with the character $x\to e^{ixt}$.  In that case
\[\hat{f}(t)=\int f(x)e^{-ixt}\,dx,\hspace{.2in} 
\check{\mu}(x)=\int e^{itx}\,d\mu(t).   \]
This is just the usual Fourier transform.  

Now $M(\widehat{G})$ is a convolution Banach algebra, and contains $L^{1}(\widehat{G})$ as a closed ideal.  The {\em Fourier-Stieltjes algebra} $B(G)$ is defined to be 
$M(\widehat{G})\spcheck$ while the {\em Fourier algebra}, $A(G)$, is defined to be $L^{1}(\widehat{G})\spcheck$.  These are algebras of continuous bounded functions on $G$ under pointwise product, with $A(G)\subset C_{0}(G)$ an ideal in $B(G)$.  Then $B(G)$ is a Banach algebra under the $M(\widehat{G})$ norm, and $A(G)\cong L^{1}(\widehat{G})$ is a closed ideal in $B(G)$.  There is a substantial literature on $A(G)$ in the abelian case - see for example \cite{Reiter}.  We note some important properties of $A(G)$.  Firstly, the set of functions of the form $f*g$, where $f,g\in C_{c}(G)$, is dense in $A(G)$.  Next, $A(G)$ is a dense subalgebra of $C_{0}(G)$ under the $\sup$-norm.  Further, $\norm{.}_{\infty}\leq \norm{.}$ on $A(G)$.
Then the Gelfand space (character space) of $A(G)=L^{1}(\widehat{G})$ is $G$, the elements of $G$ acting by point evaluation.  Lastly, $A(G)$ has a bounded approximate identity. 

To extend the notions of $A(G), B(G)$ to the non-abelian case, where duality theory is more complicated, we need to be able to interpret the norms on these algebras without reference to a dual space $\widehat{G}$.  This was done by Eymard (\cite{Eymard}) and will now be described.  It provides an excellent example of {\em quantization}, a process in which one procedes from the commutative situation to the non-commutative by replacing functions by Hilbert space operators. The following four spaces in the commutative case will be replaced by their natural ``non-commutative'' versions: 
\[ A(G)=L^{1}(\widehat{G}), B(G)=M(\widehat{G}), C_{0}(\widehat{G}), L^{\infty}(\widehat{G}).        \]
These versions will involve $G$ only, so that no duality comes in explicitly, and they make sense in the non-commutative case.  

For general $G$, we define the left regular representation $\pi_{2}$ of $G$ on $L^{2}(G)$ by:
$\pi_{2}(x)f(t)=f(x^{-1}t)$.  There is also the universal representation $\pi_{univ}$ of $G$ on a Hilbert space $H_{univ}$.
Every unitary representation of $G$ determines by integration a non-degenerate $^{*}$-representation of $C_{c}(G)$.  The norm closure of $\pi_{2}(C_{c}(G))$ is called the reduced $\cs$ $C_{red}^{*}(G)$ of $G$, while that of $\pi_{univ}(C_{c}(G))$ is called the {\em universal} $\cs$ of $G$.  The von Neumann algebra generated by $C_{red}^{*}(G)\subset B(L^{2}(G))$ is denoted by $VN(G)$.  

Suppose now that $G$ is abelian.  The key to the desired quantization is {\em Plancherel's theorem}: {\bf $f\to \check{f}$ is an isometry from 
$L^{2}(G)$ onto $L^{2}(\widehat{G})$.}  This induces an isomorphism $\Phi$ from $B(L^{2}(\widehat{G}))$
onto $B(L^{2}(G))$. Regard $C_{0}(\widehat{G}), L^{\infty}(\widehat{G})$ as $C^{*}$-subalgebras of $B(L^{2}(\widehat{G}))$ by having them act as multiplication operators on $L^{2}(\widehat{G})$.  Then $\Phi$ identifies $C_{0}(\widehat{G})$ with $C_{red}^{*}(G)$ and 
$L^{\infty}(\widehat{G})$ with $VN(G)$.  The predual $A(G)=L^{1}(\widehat{G})$ of $L^{\infty}(\widehat{G})$ goes over to the predual $VN(G)_{*}$ of $VN(G)$.  For general $G$, we then define $A(G)$ to be $VN(G)_{*}$.  To define $B(G)$, in the abelian case, it is obvious that $C_{0}(\widehat{G})=C^{*}(G)$, and so $B(G)=M(\widehat{G})=C_{0}(\widehat{G})^{*}$ is just the Banach space dual $C^{*}(G)^{*}$.  For general $G$ we define $B(G)$ to be $C^{*}(G)^{*}$.  In the natural way, $A(G)$ is a subspace of $B(G)$.  In fact since $C_{red}^{*}(G)$ is a homomorphic image of $C^{*}(G)$, it follows that $A(G)$ is closed in $B(G)$.

The space $B(G)$ can be regarded as a subspace of $C(G)$ as follows. Every $\phi\in B(G)$ is {\em coefficient} of $\pi_{univ}$, i.e. there exist vectors $\xi,\eta\in H_{univ}$ such that for all 
$f\in C_{c}(G)$, we have $\phi(f)=\lan \pi_{univ}(f)\xi,\eta\ran$.  We will write (following J. Renault (\cite{Ren})) $\phi=(\xi,\eta)$. The function $\phi$ is called {\em positive definite} if we can take $\xi=\eta$.  This property can be abstractly characterized as follows: 
\begin{equation}
\iint \phi(y^{-1}x)f(y)\ov{f(x)}\,d\la(x)\,d\la(y)\geq 0.  \label{eq:posd}
\end{equation}
for all $f\in C_{c}(G)$.  

An important subset of $B(G)$ is the set of {\em positive definite functions} on $G$.  Every element of $B(G)$ is (by polarization) a linear combination of positive definite functions.  In the abelian case, $P(G)$ is identified with the set of positive measures in $M(\widehat{G})$.

The $C^{*}(G)^{*}$ norm on $B(G)$ can (\cite{Eymard}) be usefully expressed in terms of coefficients: $\norm{\phi}$ is just $\inf\norm{\xi}\norm{\eta}$, the $\inf$ being taken over all pairs $\xi,\eta$ for which $\phi=(\xi,\eta)$. (There is another possible way of defining the norm 
of $B(G)$ that will be considered later in the groupoid context.)
Using direct sums and tensor products of representations of $G$, one obtains that $B(G)$ is a commutative Banach algebra.  

Turning to $A(G)$, it can be defined as the closure of the subspace of $B(G)$ spanned by the coefficients of $\pi_{2}$.  In fact, Eymard shows that $A(G)$ is exactly the set of functions
$(f,g)=\ov{g}*f^{\dagger}$ where $f,g\in L^{2}(G)$.  (Here, for $f\in L^{2}(G)$, $f^{\dagger}(x)=f(x^{-1}), f^{*}(x)=\ov{f(x^{-1})}$.)
An important result of Godement (\cite{Godement},\cite[Theorem 13.8.6]{D2}) says that if $\phi\in L^{2}(G)\cap P(G)$, then $\phi=f*f^{*}$ for some $f\in L^{2}(G)$.  It follows that $B(G)\cap C_{c}(G)$ is a dense subspace of $A(G)$, and $A(G)$ is a closed ideal in $B(G)$.  In particular, $A(G)$ is a commutative Banach algebra.  Eymard proves (\cite{Eymard})
his remarkable duality theorem for a general locally compact group, viz.  {\em the character space of $A(G)$ is identified with $G$, the elements of $G$ acting as characters by point evaluation}.  Walter (\cite{Wal}) showed that both $A(G)$ and $B(G)$ as Banach algebras determine the group $G$.  The important properties of $A(G)$, given in the paragraph above defining $B(G)$ in the abelian case, all hold except that $A(G)$ need not have a bounded approximate identity.  (We will return to this below in our discussion of {\em amenability}.)

Since the work of Eymard, substantial progress has been made in developing the theory of $B(G)$, $A(G)$, and before leaving the group situation and turning to locally compact groupoids, we will very briefly describe some of the main themes of this progress and cite some of the papers involved.  I am grateful to Brian Forrest for helpful information about the literature on $A(G)$. The account given below is not comprehensive, and the present writer apologizes to
all authors whose worthy contributions have been omitted.  For useful surveys of the field, the reader is referred to \cite{Pier} (for work up to 1984), and \cite{Lau94}. 

A striking feature of the theory is the relation that the amenability of $G$ has to properties of $A(G)$.  Recall that a locally compact group $G$ is called {\em amenable} if there exists an invariant mean on $L^{\infty}(G)$.  There are many other characterizations of amenable groups, and the rich phenomenon of amenability is discussed in detail in \cite{Pier,Pat1}.  
Leptin (\cite{Leptin}) showed that the Banach algebra $A(G)$ has a bounded approximate identity if and only if $G$ is amenable.  Further (\cite{Cowling,Herz2,McKennon,Renaud,Losert}) the multiplier algebra of $A(G)$ is canonically isomorphic to $B(G)$ if and only if $G$ is amenable.
There are also $L_{p}-$ versions of $A(G)$, $B(G)$ due to Herz (\cite{Herz}) which have been further investigated (e.g. in \cite{Pier,Granirer85,Granirer87}).

Banach algebra amenability of Fourier algebras has proved intriguing.  B. E. Johnson (\cite{John}) showed that the amenability of $G$ does not entail the amenability of the Banach algebra $A(G)$.  (In fact, he showed that it fails for the compact group $SU(2,\C)$.)  Now $A(G)$ is more than a Banach algebra - it is a completely contractive Banach algebra, and a theory of operator amenability can be developed in a natural way. Indeed, while amenability for a Banach algebra $A$ means (\cite{Johnmem}) that every bounded derivation from $A$ into a dual Banach $A$-module is inner, so operator amenability for a completely contractive Banach algebra $A$ means that every completely bounded derivation from $A$ into the dual of an operator $A$-bimodule is inner.  Remarkably, Ruan (\cite{Ruan95}) showed that $G$ is operator amenable if and only if its Fourier algebra $A(G)$ is operator amenable.  See the book by Effros and Ruan for a comprehensive treatment of this and related theorems (\cite[Ch. 16]{Effrosruan}).  The theory has been extended to Kac algebras (\cite{RuanX,Krausr}).  We note that in the operator space context, $A(G)$ is regarded as the ``convolution algebra of the dual quantum group'' (\cite{Effrosr}). 

On the Banach-Hochschild cohomology of $A(G)$, the reader is referred to the papers \cite{Laul,Laulw,Lauloy,Forrest88,Forrest92a}, and for information on topological centers, to \cite{Bakerlp,Lauloy}.  (The topological center of a Banach algebra $A$ is the set of weak$^{*}$-bicontinuous elements of $A^{**}$.) The dual and second dual spaces of $A(G)$ are studied in 
\cite{Granirer87,Granirer96b,Granirer97a,Granirer97b,Granirer97c,Forrest91,Forrest93,Forrest97}.
On the ideal structure of Fourier algebras, see \cite{Forrest90,Forrest92b,Woods}.  For other studies of $A(G)$, see \cite{Forrest94a,Forrest98,Granirer96a,Hu1,Hu2,Hu3,Hu4,Hu5,Kaniuthl0,Kaniuthl1,Lauu}.

\section{Locally compact groupoids}

Accounts of the theory of locally compact groupoids are given in the books of Jean Renault (\cite{rg}) and the present writer (\cite{Pat1}).  (See also the CBMS conference lectures by Paul Muhly (\cite{MuhlyTCU}).)  We summarize here the basic theory that will be needed in our discussion of $B(G)$, $A(G)$ in the groupoid case.
(There is also an important theory of Lie groupoids - this is discussed in \cite{Lands,Mackenzie,Pat1}.)

A {\em groupoid} is most simply defined as a small category with inverses.
Spelled out axiomatically, a groupoid
is a set $G$ together with a subset $G^{2} \subset G\times G$, a ``product map''
$m:G^{2}\to G$, where we write $m(a,b)=ab$,
and an inverse map $i:G\to G$, where we write $i(a)=a^{-1}$
and  where $(a^{-1})^{-1}=a$ for all $a\in G$,  such that:
\begin{enumerate}
\item if $(a,b), (b,c)\in G^{2}$, then $(ab,c), (a,bc)\in G^{2}$
and
\[     (ab)c=a(bc);        \]
\item $(b,b^{-1})\in G^{2}$ for all $b\in G$, and if $(a,b)$
belongs to $G^{2}$, then
\[     a^{-1}(ab) = b \hspace{.2in} (ab)b^{-1} = a. \]
\end{enumerate}

We define the {\em range} and {\em source} maps $r:G\to G^{0}$, $s:G\to G^{0}$
by setting $r(x)=xx^{-1}, s(x)=x^{-1}x$.
The {\em unit space} $G^{0}$ is defined to be $r(G)$ $(=s(G))$, or equivalently,
the set of idempotents $u$ in $G$. Each of the maps $r,s$ fibers the groupoid $G$ over
$G^{0}$ with fibers $\{G^{u}\}, \{G_{u}\}$ (where $u\in G^{0}$), so that $G^{u}=r^{-1}(\{u\})$ and
$G_{u}=s^{-1}(\{u\})$. Note that $(x,y)\in G^{2}$ if and only if $s(x)=r(y)$.

Intuitively, then, a groupoid is a set with a partially defined product and with inverses so that the usual group axioms hold whenever they make sense.  Groupoids give the algebra of local symmetry while groups give that of global symmetry. (See, for example, \cite{Weinstein}). Groupoids are becoming more and more prominent in analysis.  Connes's non-commutative geometry (\cite{Connesbook}) makes extensive use of them and this has been an important motivation for their study.  Examples of groupoids are:
\bi                                                  
\item[(a)]locally compact groups
\item[(b)]equivalence relations
\item[(c)]tangent bundles
\item[(d)]the tangent groupoid (e.g. \cite{Connesbook})
\item[(e)]holonomy groupoids for foliations (e.g. \cite{Connesbook})
\item[(f)]Poisson groupoids (e.g. \cite{Weinp})
\item[(g)]transformation groups (e.g. \cite{MuhlyTCU})
\item[(h)] graph groupoids (e.g. \cite{KumRaeRen,Pat4})
\ei

As a simple, helpful example of a groupoid, consider (b) above.
Let $R$ be an equivalence relation on a set $X$.  Then $R$ is a groupoid under the following operations:
$(x,y)(y,z)=(x,z), (x,y)^{-1}=(y,x)$. Here, $G^{0}=X$ (=diagonal of $X\x X$) and 
$r((x,y))=x, s((x,y))=y$.  So $R^{2}=\{((x,y),(y,z)): (x,y), (y,z)\in R\}$.  When $R=X\x X$, then $R$ is called a {\em trivial} groupoid.

A special case of a trivial groupoid is $R=R_{n}=\{1,2,\ldots ,n\}\x \{1,2,\ldots ,n\}$.  (So every $i$ is equivalent to every $j$.)  Identify $(i,j)\in R_{n}$ with the matrix unit $e_{ij}$.  Then the groupoid $R_{n}$ is just matrix multiplication except that we only multiply $e_{ij}, e_{kl}$ when $k=j$, and $(e_{ij})^{-1}=e_{ji}$.  We do not really lose anything by restricting the multiplication, since the pairs $e_{ij}, e_{kl}$ excluded from groupoid multiplication just give the 0 product in normal algebra anyway.

For a groupoid $G$ to be a {\em locally compact} groupoid means what one would expect.  The groupoid $G$ is required\footnote{In a number of contexts, e.g. that involving the holonomy groupoids for foliations, the Hausdorff condition is too strong and locally Hausdorff groupoids need to be considered.  For a discussion of the theory of such groupoids, the reader is referred to \cite{Pat2,Khoshskand}.  Only locally compact Hausdorff groupoids will be considered in the present paper.} to be a (second countable) locally compact Hausdorff space, and the product and inversion maps are required to be continuous.  Each $G^{u}$ as well as the unit space $G^{0}$ is closed in $G$.
What replaces left Haar measure on $G$ is a system of measures $\la^{u}$ $(u\in G^{0})$, where $\la^{u}$ is a positive regular Borel measure on $G^{u}$ with dense support.  In addition, the $\la^{u}$'s  are required to vary continuously (when integrated against $f\in C_{c}(G)$)
and to form an {\em invariant} family in the sense that for each $x$, the map
$y\to xy$ is a measure preserving homeomorphism from $G^{s(x)}$ onto $G^{r(x)}$. Such a system $\{\la^{u}\}$ is called a {\em left Haar system} for $G$.  Not all locally compact groupoids possess a left Haar system, and even if there exists a left Haar system, it need not be unique.  However, most of the locally compact groupoids that arise in practice have natural left Haar systems.  All locally compact groupoids in the present survey are assumed to have a left Haar system.The presence of a left Haar system on $G$ has topological implications: it implies that the range map $r:G\to G^{0}$ is open.

For such a $G$, the vector space $C_{c}(G)$ is a convolution $^{*}$-algebra, where for $f,g\in C_{c}(G)$:
\[  f*g(x)=\int f(t)g(t^{-1}x)\, d\la^{r(x)}(t),\hspace{.2in}  f^{*}(x)=\ov{f(x^{-1})}.       \]
We take $C^{*}(G)$ to be the enveloping $C^{*}$-algebra of $C_{c}(G)$ (representations are required to be continuous in the inductive limit topology).  Equivalently, it is the completion of $\pi_{univ}(C_{c}(G))$ where $\pi_{univ}$ is the universal representation of $G$.  For example, if $G=R_{n}$, then $C^{*}(G)$ is easy to guess! - it is just the finite dimensional algebra $C_{c}(G)=M_{n}$, the span of the $e_{ij}$'s. 

The class of locally compact groupoids which corresponds to discrete groups in the group category is that of the {\em r-discrete} groupoids.  We can think of such a groupoid $G$ as one for which the unit space $G^{0}$ is an arbitrary locally compact space and the fibers $G^{u}$ are discrete with local sections.  The latter condition is more precisely expressed as saying that every $x\in G$ has an open neighborhood $U$ for which $r(U)$ is open in $G$ and the map $r_{\mid U}$ is a homeomorphism from $U$ onto $r(U)$.  (The unit space $G^{0}$ is therefore {\em open} in $G$.) The canonical left Haar system for $G$ is that for which each $\la^{u}$ is counting measure on $G^{u}$.  Of course, every discrete groupoid is r-discrete.  Other examples are transformation groupoids for which the acting group is discrete, and certain kinds of holonomy groupoids. 

We will require in \S 4 the notion of a {\em bisection} $a$ of a locally compact groupoid $G$.  Here $a$ is a pair of homeomorphisms $u\to a^{u}\in G^{u}$, $u\to a_{u}\in G_{u}$ from $G^{0}$ onto a subset $A$ of $G$.  So $A$ is simultaneously an r-section and an s-section and determines and is determined by the bisection $a$.  The set of bisections of $G$ forms, by multiplying and inverting the $A$'s, a group, and is denoted by $\Ga$.  If $G$ is a locally compact group, then trivially $\Ga=G$.   

As in the case of a locally compact group, there is a reduced $\cs$
$C_{red}^{*}(G)$ for $G$ which is defined as follows.  For each $u\in
G^{0}$, we first define
a representation $\pi_{u}$ of $C_{c}(G)$
on the Hilbert space $L^{2}(G,\la_{u})$.  To this end,
regard $C_{c}(G)$ as a dense subspace of $L^{2}(G,\la_{u})$ and define
for $f\in C_{c}(G), \xi\in C_{c}(G)$,
\begin{equation}
     \pi_{u}(f)(\xi)=f*\xi\in C_{c}(G).       \label{eq:red}
\end{equation}
Then $\pi_{u}(f)$ extends to a bounded linear operator on $L^{2}(G,\la_{u})$.
The reduced $\cs$-norm on $C_{c}(G)$ is then (e.g. \cite[p.108]{Pat2})
defined by:
\[        \norm{f}_{red}=\sup_{u\in G^{0}}\norm{\pi_{u}(f)}       \]
and $C_{red}^{*}(G)$ is defined to be the completion of $C_{c}(G)$ under this norm.

We next recall some details concerning the disintegration of representations of $C_{c}(G)$.  The theorem is due to J. Renault (\cite{Renrep}).  A detailed account of the theorem is given in the book of Paul Muhly (\cite{MuhlyTCU}).  Let $\Phi$ be a representation of $C_{c}(G)$ on a Hilbert space $H$ which is continuous in the inductive limit topology.  Then $\Phi$ disintegrates as follows.  

There is a probability measure $\mu$ on $G^{0}$ which is {\em quasi-invariant} (in the sense defined below).  Associated with $\mu$ is a positive regular Borel measure $\nu$ on $G$ defined by: $\nu=\int\la^{u}\,d\mu(u)$.
The measure $\nu^{-1}$ is the image under $\nu$ by inversion: precisely, $\nu^{-1}(E)=\nu(E^{-1})$.  There is also a measure $\nu^{2}$ on $G^{2}$ given by: $\nu^{2}=\iint \la_{u}\x \la^{u}\,d\mu(u)$. The quasi-invariance of $\mu$ just means that $\nu$ is equivalent to $\nu^{-1}$.  The {\em modular function} $D$ is defined as the Radon-Nikodym derivative $d\nu/d\nu^{-1}$. The function $D$ can be taken to be Borel (\cite{Hahn,RamsayJFA,MuhlyTCU}), and satisfies the 
properties: $D(x^{-1})=D(x)^{-1}$ $\nu$-almost everywhere, and 
$D(xy)=D(x)D(y)$ $\nu^{2}$-almost everywhere.  Let $\nu_{0}$ be the measure on $G$ given by: $d\nu_{0}=D^{-1/2}d\nu$.

Next, there exists (e.g. \cite[p.91]{Pat2}) a measurable Hilbert bundle $(G^{0},\mfK,\mu)$ 
($\mfK=\{\mfK^{u}\}_{u\in G^{0}}$) and a {\em $G$-representation} $L$ on $\mfK$.  This means that each $L(x)$ ($x\in G$) is a linear isometry from $\mfK^{s(x)}$ onto $\mfK^{r(x)}$ and $L(x)$ is the identity map if $x\in G^{0}$.  Further, the map $x\to L(x)$ is multiplicative $\nu^{2}$-almost everywhere and inverse preserving 
$\nu$-almost everywhere.  Lastly, for every pair $\xi,\eta$ of square integrable sections of $\mfK$, it is required that the function $x\to (L(x)\xi(s(x)),\eta(r(x)))$ be $\nu$-measurable.  The representation $\Phi$ of $C_{c}(G)$ is then given by:
\begin{equation}  
\lan \Phi(f)\xi,\eta\ran = \int f(x)(L(x)\xi(s(x)),\eta(r(x)))\,d\nu_{0}(x).   \label{eq:disinteg}
\end{equation}
We will refer to the triple $(\mu,\mfK,L)$ as a {\em measurable $G$-Hilbert bundle}.
The notion of a {\em continuous} $G$-Hilbert bundle, that will be needed in \S 5, is defined in the obvious way. 

Significant progress has been made in recent years in developing theories of the Fourier-Stieltjes and Fourier algebras for locally compact groupoids, and these theories will be described in the rest of the paper.  The main papers involved are \cite{Ren,Ramwal,Pat3}, and these will be considered in turn.  As we will see, the theories relate to one another.

\section{The Fourier-Stieltjes and Fourier algebras for a measured groupoid}

This section discusses the paper \cite{Ren} of Jean Renault that deals with the Fourier algebra of a {\em measured groupoid}.  We can take the latter to mean (\cite{RamsayJFA})
that $G$ is a locally compact groupoid with left Haar system $\{\la^{u}\}$ and fixed quasi-invariant measure $\mu$ with modular function $D$.  The theory is a far reaching extension of the discussion of these algebras in \S 1 for locally compact groups, and as in the group case, has an operator algebraic character.  (In the group case, $G=G^{e}$ where $e$ is the identity of $G$, $\la^{e}$ is a left Haar measure and $\mu$ is the point mass at $e$.)  

As we saw in \S 1, the Fourier-Stieltjes algebra $B(G)$ of a locally compact group is the space of coefficients $(\xi,\eta)$ of Hilbert space representations of $G$.  In the groupoid case, the Fourier-Stieltjes algebra $B_{\mu}(G)$ is defined to be the space of coefficients $\phi=(\xi,\eta)$ where $\xi,\eta$ are $L^{\infty}$-sections for some measurable $G$-Hilbert bundle $(\mu,\mfK,L)$.  So for $x\in G$, 
\[           \phi(x)=(L(x)\xi(s(x)),\eta(r(x))).    \]  
Clearly, $\phi$ belongs to $L^{\infty}(G)=L^{\infty}(G,\nu)$.

As in the group case, the set $P_{\mu}(G)$ of positive definite functions in $L^{\infty}(G)$ plays an important role.  A function $\phi\in L^{\infty}(G)$ is called {\em positive definite} if and only if for all $u\in G^{0}$,
\[ \iint \phi(y^{-1}x)f(y)\ov{f(x)}\,d\la^{u}(x)\,d\la^{u}(y)\geq 0.  \label{eq:posd2}   \]
An argument of Ramsay and Walter (\cite{Ramwal}) gives that $\phi$ is positive definite if and only if $\phi$ is of the form $(\xi,\xi)$ for some $\xi$.  By polarization, $B_{\mu}(G)$ is the span of $P_{\mu}(G)$.

Also, as in the group case, the norm $\norm{\phi}$ of $\phi\in B_{\mu}(G)$ is defined to be $\inf\norm{\xi}\norm{\eta}$ over all representations $\phi=(\xi,\eta)$.  Using a groupoid version of Paulsen's ``off-diagonalization technique'' (\cite[Th. 7.3]{Paulsen}, \cite[Th. 5.3.2]{Effrosruan}), Renault shows that $B_{\mu}(G)$ is a commutative Banach algebra.

The Fourier algebra $A_{\mu}(G)$ is defined to be the closed linear span in $B_{\mu}(G)$ of the coefficients of the {\em regular} representation of $G$ on the $G$-Hilbert bundle 
$\{L^{2}(G^{u})\}_{u\in G^{0}}$, the $G$-action being given by left translation. (We will meet the continuous version of this Hilbert bundle again in \S 5.) Renault shows that $A_{\mu}(G)$ is a closed ideal in $B_{\mu}(G)$.  The {\em Hilbertian} functions of Grothendieck (\cite{Groth}) give an example of a $B_{\mu}(G)$ (with $G$ a trivial groupoid).  

Renault then goes on to examine duality theory using operator algebra and operator space techniques.  Corresponding to the universal $\cs$ $C^{*}(G)$ in the group case is the universal $C_{\mu}^{*}(G)$ in the measured groupoid case.  The latter is the completion of $C_{c}(G)$ under the largest $C^{*}$-norm coming from some measurable $G$-Hilbert bundle $(\mu,\mfK,L)$.  The reduced $\cs$ $C_{red}^{*}(G)$ and the von Neumann algebra $VN(G)$ (both depending on $\mu$) are defined using the regular representation in the natural way.   

In the group case, we saw that $B(G)=C^{*}(G)$.  In the measured groupoid case, the situation is more subtle.  Some operator space notions are required.  Recall that any Hilbert space $H$ can be regarded as a operator space by identifying it in the obvious way with a subspace of $\mfB(\C,H)$: each $\xi\in H$ is identified with the map $a\to a\xi$ ($a\in \C$).  Further, trivially, $H^{*}$ is an operator space as a subspace of $\mfB(H,\C)$.  In the measured groupoid situation, Renault shows that the operator spaces $L^{2}(G^{0})$ and $C_{\mu}^{*}(G)$ are completely contractive left $L^{\infty}(G^{0})$ modules, and $L^{2}(G^{0})^{*}$ is a completely contractive right $L^{\infty}(G^{0})$ module.  If $E$ is a right and $F$ is a left $A$-operator module ($A$ a $\cs$) then the Haagerup tensor norm is determined on the algebraic tensor product
$E\odot_{A} F$ by setting
$\norm{u}=\sum_{i=1}^{n}\norm{e_{i}}\norm{f_{i}}$ over all representations 
$u=\sum_{i=1}^{n}e_{i}\otimes_{A} f_{i}$.  The completion
(\cite[Ch. 2]{BMP}) $E\otimes_{A} F$ of $E\odot_{A} F$ is called
the {\em module Haagerup tensor product} of $E$ and $F$ over $A$.

One then forms the module Haagerup tensor product
$X(G)=L^{2}(G^{0})^{*}\otimes C_{\mu}^{*}(G)\otimes L^{2}(G^{0})$ over $L^{\infty}(G^{0})$.  Renault proves that $X(G)^{*}=B_{\mu}(G)$.  Under this identification, each 
$\phi=(\xi,\eta)$ goes over to the linear functional $a^{*}\otimes f\otimes b\to
\int \ov{a\circ r}(\phi f)b\circ s\,d\nu$ ($f\in C_{c}(G)$).  

Renault also gives a characterization of $A_{\mu}(G)$ which is the groupoid version of the group result (\S 1) that $A(G)=VN(G)_{*}$.  In fact the predual $VN(G)_{*}$ of $VN(G)$ is completely isometric to the module Haagerup tensor product $L^{2}(G^{0})^{*}\otimes A_{\mu}(G)\otimes L^{2}(G^{0})$ over $L^{\infty}(G^{0})$.  He also shows that the analogue of the group multiplier result, described in \S 1, holds, viz. {\em if the measure groupoid $G$ is amenable (which can be defined as saying that the trivial representation is weakly contained in the regular representation) then $B_{\mu}(G)$ is the multiplier algebra of $A_{\mu}(G)$.}  He then uses these results on $A(G), B(G)$ to investigate absolute Fourier multipliers for r-discrete groupoids, generalizing results of Pisier and Varopoulos. 

Leptin's result (\S 1) relating the amenability of a locally compact group to the existence 
of a bounded approximate identity in the Fourier algebra has been generalized to the measured groupoid context by Jean-Michel Vallin (\cite{Vall1,Vall2}) using Hopf-von Neumann bimodule structures.

\section{The Borel Fourier-Stieltjes algebra}

In this section we describe some of the results of the paper \cite{Ramwal} by Ramsay and Walter.  The setting for the paper is that of a locally compact groupoid $G$.  The groupoid is not considered to be a measured groupoid, so that a quasi-invariant measure $\mu$ on $G^{0}$ is not specified in advance.   We will also consider briefly a paper by K. Oty at the end of the section.

Since in the group case, $B(G)$ is the span of the set $P(G)$ of continuous positive definite functions, it would be natural to use this as the definition of $B(G)$ in the groupoid case.  We will in fact effectively use this definition of $B(G)$ in \S 5. However, there are examples (\cite[\S 7]{Ramwal}) in which the span of the continuous positive definite functions on $G$ is not complete in the norm (defined below) on $B(G)$.  Instead, Ramsay and Walter consider the span $\mathcal{B}(G)$ of the set $\mathcal{P}(G)$ of bounded {\em Borel} positive definite functions on $G$.  There are examples of {\em continuous} functions in $\mathcal{B}(G)$ which cannot be expressed as linear combination of elements of $P(G)$.  In $\mathcal{B}(G)$, we identify two functions that are 
$\la^{\mu}$-equal for every quasi-invariant measure $\mu$ on $G^{0}$.

The norm on $\mathcal{B}(G)$ is defined in terms of a remarkable {\em completely bounded multiplier norm} on the $\cs$ $M^{*}(G)$. (This norm does not seem to have been considered in 
theory of the Fourier-Stieltjes algebra in the group case.)
In more detail, the universal representation $\pi_{univ}$ extends canonically to a representation of the convolution algebra $B_{c}(G)$ of compactly supported, bounded Borel functions on $G$.  Let $M^{*}(G)$ be the completion of $\pi_{univ}(B_{c}(G))$.  Then $M^{*}(G)$ is a $\cs$.  Each $\phi\in \mathcal{B}(G)$ acts as a 
multiplier $T_{\phi}$ on $M^{*}(G)$ by extending its action from $\pi_{univ}(B_{c}(G))$ by continuity.  Here, for $f\in B_{c}(G)$, 
$T_{\phi}\pi_{univ}(f)=\pi_{univ}(\phi f)$ ($\phi f$ pointwise multiplication on $G$).  Ramsay and Walter show that $T_{\phi}$ is a completely bounded operator on $M^{*}(G)$.  They define $\norm{\phi}_{cb}=\norm{T}_{cb}$, and show that $\mathcal{B}(G)$ is a Banach algebra under $\norm{.}_{cb}$.  (Below, we will identify the version of $\norm{.}_{cb}$ for a  locally compact abelian group in the continuous context.)

Ramsay and Walter give a partial analogue of the locally compact groupoid result (\S 1) that 
$B(G)=C^{*}(G)^{*}$.  Let $X$ be the one point compactification of $G^{0}$.  They construct two $\css$ $C^{*}(G,X)$, $M^{*}(R,X)$
which are $C(X)$-submodules.  Here $R$ is the orbit equivalence relation on $G^{0}$: so $u\sim v$ if and only if there exists $x\in G$ such that $r(x)=u, s(x)=v$.  They show using a technical argument that each $\phi\in \mathcal{B}(G)$ gives a completely bounded $C(X)$-bimodule map 
$S_{\phi}: C^{*}(G,X)\to  M^{*}(R,X)$.  (When $G$ is a group, then $S_{\phi}$ is just the canonical linear functional $\phi: C^{*}(G)\to \C$.)  It would be interesting to know how this result relates to the corresponding duality result involving $X(G)$ for measured groupoids in \S 3.

Karla Oty (\cite{Oty}) investigates the space $\mathcal{B}(G)\cap C(G)$ and the set
$P(G)\subset \mathcal{P}(G)$ of continuous positive definite functions.  She shows among other results that $P(G)$ separates the points of $G$, and that if $G$ is r-discrete then we do not have to work with equivalence classes in $\mathcal{B}(G)$.  Further, if the set of cardinalities of the 
$G^{u}$'s ($u\in G^{0}$) is bounded by a positive integer, then $\mathcal{B}(G)\cap C(G)=C(G)$.

\section{The continuous Fourier-Stieltjes and Fourier algebras} 

In this section, we consider some of the results in the paper \cite{Pat3} by the present writer.  
In contrast to the measurable theory of Renault, considered in \S 2, and the Borel theory of Ramsay and Walter, considered in \S 3, the theory of this section is a continuous one.  For motivation, in \S 3, in the context of a measured groupoid, $A_{\mu}(G)$ was defined as the closed linear span in $B_{\mu}(G)$ of the coefficients of the regular $G$-Hilbert bundle $\{L^{2}(G^{u})\}$.  Now this is a {\em continuous} Hilbert bundle, the continuous sections being determined in the obvious way by $C_{c}(G)$.  So
it is natural to consider the coefficients of {\em continuous} $G$-Hilbert bundles as our space $B(G)$.  The elements of $B(G)$ are then all continuous as in the group case. 
As we will see, the ideas and techniques from \cite{Ren,Ramwal} play a fundamental role in the theory.  We will have occasion to consider Hilbert modules, and so in accordance with the usual practice (e.g. \cite{Lance}), for the rest of the paper, Hilbert spaces and modules will be taken to be conjugate linear in the first variable.

We now define the Fourier-Stieltjes algebra $B(G)$ in more detail.  Given a continuous $G$-Hilbert bundle $\mfH$, we consider the Banach space $\De_{b}$ of continuous, bounded sections of $\mfH$. 
For $\xi,\eta\in \De_{b}$, define (as in \S 3) the {\em coefficient} $(\xi,\eta)\in C(G)$ by: $(\xi,\eta)(u)=(L_{x}\xi(s(x)),\eta(r(x)))$ where $x\to L_{x}$ is the $G$-action on $\mfH$.  Then $B(G)$ is defined to be the set of all such coefficients, coming from all possible continuous $G$-Hilbert bundles.  As in \cite{Ren}, $B(G)$ is an algebra over $\C$ and the norm of $\phi\in B(G)$ is defined to be $\inf\norm{\xi}\norm{\eta}$, the $\inf$ being taken over all representations  
$\phi=(\xi,\eta)$. Then $B(G)\subset C(G)$, and $\norm{.}_{\infty}\leq \norm{.}$.  An argument similar to that of \cite{Ren} shows that $B(G)$ is a commutative Banach algebra.

From an argument of \cite{Ramwal}, $P(G)\subset B(G)$.
Clearly, using the polarization for sections of continuous Hilbert bundles, every $\phi\in B(G)$ is a linear combination of continuous positive definite elements of $B(G)$ and conversely.  From the discussion in \S 4, we see that $\mathcal{B}(G)\cap C(G)\neq B(G)$ in general.  Intuitively, $\mathcal{B}(G)$ captures the representation theory of $G$ but we need to consider Borel functions.  On the other hand, $B(G)$ may not capture all of the representation theory of $G$ (though it does capture at least the regular representation) but has the advantage that the functions involved, are, as in the group case, continuous.  (Our $B(G)$ is the same as Oty's $B_{1}(G)$ (\cite{Oty}).)  

We defined the norm on $B(G)$ using Hilbert bundles in a way similar to that in which Renault defined the norm on $B_{\mu}(G)$, and it is natural to enquire if the method of Ramsay and Walter for norming $\mathcal{B}(G)$ can also be applied in the present context.  This is in fact the case.  The difference is that we regard $\phi\in B(G)$ as a multiplier on $C_{c}(G)$, this action extending to a completely bounded map $T_{\phi}$ on $C^{*}(G)$ instead of on $M^{*}(G)$.  
We can then define a norm $\norm{.}_{cb}$ on $B(G)$ by taking $\norm{\phi}_{cb}=\norm{T_{\phi}}_{cb}$.  Following along the same lines as the corresponding argument in \cite{Ramwal}, one can show that $B(G)$ is a Banach algebra under $\norm{.}_{cb}$.  Since it is not difficult to show that $\norm{.}_{cb}\leq \norm{.}$, it follows by Banach's isomorphism theorem that the two norms are equivalent on $B(G)$.  It is an open question if the two norms actually coincide in general.  It is shown, using a result of Paulsen on Schur multipliers (\cite{Paulsen}) that they do coincide when $G=R_{n}$, the trivial groupoid on $n$ elements (\S 2).

Let us sketch here how one can show that the two norms are the same when G is a locally compact abelian group.  As in \S 1, we can regard $\phi$ as a measure $\mu\in M(\widehat{G})$.  Recalling that $C^{*}(G)=C_{0}(\widehat{G})$, the map $T_{\phi}$ goes over to the map $T_{\mu}:C_{0}(\widehat{G})\to C_{0}(\widehat{G})$ where
\begin{equation}          \label{eq:mug}
  T_{\mu}(g)(x)=\mu*g(x)=\int g(t^{-1}x)\,d\mu(t).     
\end{equation}
Trivially, $\norm{T_{\mu}}\leq \norm{\mu}$.  With $e$ the identity of $\widehat{G}$ 
and using (\ref{eq:mug}),
\[\norm{T_{\mu}}\geq\sup_{\norm{g}=1}\mid T_{\mu}g(e)\mid =\norm{\mu}.     \]
So $\norm{T_{\phi}}=\norm{\mu}=\norm{\phi}$.  Since $C^{*}(G)$ is abelian, $\norm{\phi}_{cb}=\norm{T_{\phi}}$ (\cite[Th. 3.8]{Paulsen}). 

For the rest of the discussion, we can use either norm on $B(G)$.  We have chosen to use $\norm{.}$ since in the setting of $A(G)$ (below), the continuous $G$-Hilbert bundle 
$L^{2}(G)=\{L^{2}(G^{u})\}$ is conveniently present for norm estimates.  

Let $E^{2}$ be the space of continuous sections of $L^{2}(G)$ that vanish at $\infty$.  Let $D=C_{0}(G^{0})$.  Then $E^{2}$ 
is a Hilbert $D$-module, with inner product given by: $\lan\xi,\eta\ran(u)=(\xi(u),\eta(u))$ and left module action by: $(\xi a)(u)=a(u)\xi(u)$ ($a\in D, \xi\in E^{2}$).  (A good source of information about Hilbert $A$-modules is the book \cite{Lance} by Lance.)  

It is easy to check that if $\xi,\eta\in E^{2}$, then $\lan\xi,\eta\ran=\eta*\xi^{*}$.
For $F\in C_{c}(G)$, define $R_{F}:C_{c}(G)\to C_{c}(G)$ by right convolution:
$R_{F}f=f*F$.  Then the  map $F\to R_{F}$ is a $^{*}$-antirepresentation from $C_{c}(G)$ into 
$\mfL(E^{2})$, and the  closure of its image in $\mfL(E^{2})$ is canonically isomorphic to $C_{red}^{*}(G)$.  As in the group case, we can define $VN(G)$ to be the commutant of $C_{red}^{*}(G)$ in $\mfB(E^{2})$.  It is to be stressed that $VN(G)$ is not a von Neumann algebra in general - this is not even the case for $G=R_{n}$.  However, $VN(G)$ is always a Banach algebra, and is strongly closed in $\mfB(E^{2})$.

As in the group case, the algebra $VN(G)$ plays a useful role in the theory.  In particular, if $G$ is r-discrete then the fact that $C_{red}^{*}(G)\subset \mfL(E^{2})$ can be used to show that the natural version of Godement's theorem (\S 1) holds: i.e. that if $\phi\in C_{c}(G)\cap P(G)$, then $\phi=f*f^{*}=(f,f)$ for some $f\in E^{2}$.  (The discrete group version of this is given 
in \cite[Lemma VII.2.7]{Davidson}.)
It follows that $A(G)$ is an ideal in $B(G)$.  I do not know if these results are true for locally compact groupoids in general.

For $f,g\in C_{c}(G)$, the coefficient $(f,g)=g*f^{*}$ itself belongs to $C_{c}(G)$.  We define the Fourier algebra $A(G)$ to be {\em the closure in $B(G)$ of the subalgebra generated by the set of such elements $g*f^{*}$.}  (Another possible definition of $A(G)$ is given by Oty (\cite{Oty}).)
Of course, $A(G)$ is a commutative Banach algebra.  Further, as in the group case, $A(G)\subset C_{0}(G)$.
In the group case, $VN(G)$ is the dual of $A(G)$.  Using the $G$-Hilbert module $E^{2}$, there is a version of this for the groupoid case, though it is more involved.  For the technicalities, see \cite{Pat3}.  

The last part of \cite{Pat3} proves a duality theorem for $A(G)$ generalizing that of Eymard's in the group case.  For the latter, recall (\S 1) that the character space of $A(G)$ is just $G$, the elements of $G$ being characters under point evaluation.  I do not know if this is still true in the groupoid case.  The duality theorem of \cite{Pat3} is formulated in terms of bisections and certain multiplicative module maps from $A(G)$ into $C_{0}(G^{0})$.  In one direction, the group $\Ga$ of bisections of $G$ (\S 2) determines a pair of multiplicative maps on $A(G)$ as follows.  If $a\in \Ga$, then we define $\al^{a}:A(G)\to D$, $\bt_{a}:A(G)\to D$ by setting:
$\al^{a}(\phi)(u)=\phi(a^{u})$, $\bt_{a}(\phi)=\phi(a_{u})$.  The pair $(\al^{a},\bt_{a})$ satisfy a number of interesting properties.  For example, $\al^{a}$ is a $D$-module homomorphism for $G$, while $\bt_{a}$ is the same for $G(r)$, the groupoid $G$ with multiplication reversed.  Crucial is the fact that if $J$ is the homeomorphism $u\to v$ of $G^{0}$, where $a^{u}=a_{v}$,
then for all $\phi\in A(G)$, we have $\bt_{a}(\phi)\circ J=\al^{a}(\phi)$.  

We then consider the set $\Phi_{A(G)}$ of pairs of maps $(\al,\bt)$ satisfying these properties
abstractly.  We give $\Ga$ the topology of pointwise convergence on $G^{0}$, regarding each $a\in \Ga$ as the (single) map $u\to (a_{u},a^{u})\in G^{2}$.  Next, regarding each $(\al,\bt)\in \Phi_{A(G)}$ as the single map $(\phi,u)\to (\al(\phi)(u),\bt(\phi)(u)\in \C^{2}$ on $A(G)\x G^{0}$, we also give $\Phi_{A(G)}$ the topology of pointwise convergence.  The duality theorem of \cite{Pat3} then is that for a large class of locally compact groupoids $G$, {\em the map $a\to (\al^{a},\bt_{a})$ is a homeomorphism from $\Ga$ onto $\Phi_{A(G)}$}.  In the group case, each bisection $a$ is just a group element $x$ and the maps $\al^{a}, \bt_{a}$ coincide on $A(G)$, being just point evaluation at $x$.  Further, $\Phi_{A(G)}$ is just the space of characters of $A(G)$, and the above duality theorem reduces to Eymard's duality theorem (\S 1).

A list of open questions in the continuous theory is given in \cite{Pat3}.  This includes the natural questions as to whether $\norm{.}=\norm{.}_{cb}$ on $B(G)$ and whether the character space of $A(G)$ is identifiable with $G$.


\end{document}